\documentclass{article}

\usepackage{arxiv}

\usepackage[utf8]{inputenc} % allow utf-8 input
\usepackage[T1]{fontenc}    % use 8-bit T1 fonts
\usepackage{hyperref}       % hyperlinks
\usepackage{url}            % simple URL typesetting
\usepackage{booktabs}       % professional-quality tables
\usepackage{amsfonts}       % blackboard math symbols
\usepackage{nicefrac}       % compact symbols for 1/2, etc.
\usepackage{microtype}      % microtypography
\usepackage{lipsum}

\ifpdf \usepackage[pdftex]{graphicx} \pdfcompresslevel=9
\else \usepackage[dvips]{graphicx} \fi

\usepackage{t1enc,dfadobe}

\usepackage{egweblnk}
\usepackage{cite}

\usepackage{nicefrac}
\usepackage{amsmath}
\usepackage{bm}
\usepackage{amsfonts}
\usepackage{dsfont}
\usepackage[textsize=tiny]{todonotes}
\usepackage{subfigure}
 \usepackage{comment} 
 
 \usepackage[usenames,dvipsnames]{pstricks}
 \usepackage{epsfig}
 \usepackage{pst-grad} % For gradients
 \usepackage{pst-plot} % For axes
 \usepackage[space]{grffile} % For spaces in paths
 \usepackage{etoolbox} % For spaces in paths
 \makeatletter % For spaces in paths
 \patchcmd\Gread@eps{\@inputcheck#1 }{\@inputcheck"#1"\relax}{}{}
 \makeatother

\DeclareMathAlphabet\mathbfcal{OMS}{cmsy}{b}{n}

\title{An Introduction to the Deviatoric Tensor Decomposition in Three Dimensions
  and its Multipole Representation}

\author{
  Chiara Hergl \\
  Department of Computer Science\\
  Leipzig University \\
  Germany
  \texttt{} \\
   \And
 Thomas Nagel \\
 Geotechnical Institute\\
 Technische Universit\"at Bergakademie Freiberg\\
 Germany
  \texttt{} \\
   \And
  Gerik Scheuermann \\
  Department of Computer Science\\
  Leipzig University \\
  Germany
}

\usepackage{color}
\newcommand{\notegs}[1]{\colorbox{yellow}{[Gerik:] #1}}

\begin{document}
\maketitle

\begin{abstract}
The analysis and visualization of tensor fields is a very challenging task.
Besides the cases of zeroth- and first-order tensors, most techniques focus on symmetric second-order tensors. Only a few works concern totally symmetric tensors of higher-order. Work on other tensors of higher-order than two is exceptionally rare. We believe that one major reason for this gap is the lack of knowledge about suitable tensor decompositions for the general higher-order tensors. We focus here on three dimensions as most applications are concerned with three-dimensional space. A lot of work on symmetric second-order tensors uses the spectral decomposition. The work on totally symmetric higher-order tensors deals frequently with a decomposition based on spherical harmonics. These decompositions do not directly apply to general tensors of higher-order in three dimensions. However, another option available is the deviatoric decomposition for such tensors, splitting them into deviators. Together with the multipole representation of deviators, it allows to describe any tensor in three dimensions uniquely by a set of directions and non-negative scalars. The specific appeal of this methodology is its general applicability, opening up a potentially general route to tensor interpretation. The underlying concepts, however, are not broadly understood in the engineering community. In this article, we therefore gather information about this decomposition from a range of literature sources. The goal is to collect and prepare the material for further analysis and give other researchers the chance to work in this direction. This article wants to stimulate the use of this decomposition and the search for interpretation of this unique algebraic property. A first step in this direction is given by a detailed analysis of the multipole representation of symmetric second-order three-dimensional tensors.

\end{abstract}

% keywords can be removed

\section{Introduction}
Analysis and visualization of tensor fields focuses mainly
  on symmetric second-order fields like mechanical stress and
  strain.
There is also work on general second-order tensor fields like
  the velocity gradient in fluid flow.
But there are many more tensor fields used in natural sciences
  and engineering to describe the behaviour of matter and various fields.
However, there is only very limited work on their analysis and visualization.
In the life sciences, there is some work on the analysis and
  visualization of higher-order tensors, especially higher-order
  diffusion tensors.
However, these tensors are totally symmetric and of even-order.
If this property does not hold, literature gets exceptionally sparse.

A look at most of this work \cite{laidlaw2012new, hotz2015visualization,
  kratz2013visualization, itskov2007tensor}
%\notegs{Zitate zu den Tensor Dagstuhl Büchern und geeigneten Surveys}
  reveals that analysis and visualization of symmetric tensor fields
  nearly always uses the eigendecomposition, i.e. eigenvalues
  and eigenvectors.
For general tensors of second-order, authors split them into 
  the asymmetric part which is interpreted as vector or rotation
  and a symmetric part which again is analyzed using eigenvalues
  and eigenvectors most often.
The work on totally symmetric tensors of higher even order 
  uses a decomposition based on spherical harmonics in most cases.
Unfortunately, these decompositions do not work for general
  tensors of higher order in three dimensions.
We therefore suspect that the lack of work on the analysis and visualization of other tensors
  is caused by missing knowledge on tensor decompositions
  of these tensors.

The good news is that there is such a decomposition, namely the
  deviatoric decomposition, allowing to split each tensor
  in three dimensions into a set of deviators, i.e. traceless
  totally symmetric tensors.
Furthermore, these deviators can be uniquely represented by 
  a finite set of directions and a non-negative scalar.
The generality of this methods renders it quite appealing for study in tensor analysis and visualization.
Because this decomposition is not very well established and has gone largely unnoticed in the applied sciences, this work gives an overview of the so called deviatoric 
  decomposition.
The deviatoric decomposition is an orthogonal decomposition 
  of a tensor of arbitrary-order up to dimension three.
The deviatoric decomposition of a second-order tensor is well 
  known.
It is given by the tensor's symmetric part and a vector, representing 
  its asymmetric part.
These deviators again can be represented by the symmetrization
  of a tensor product of first-order tensors (vectors) which 
  are called multipoles.
The deviatoric decomposition and the multipole representation
  are unique, so each tensor can be represented by a unique set
  of first-order tensors, i.e. vectors. 
  
The underlying calculations go back to 
  Maxwell~\cite{Maxwell:1881} and were
summarized by Backus~\cite{backus} who also provided further 
  information and some mathematical background.
Zou et al. \cite{zou2003} gave explicit formulae to calculate
  the multipoles of a tensor.

The goal of this article is to give a summary of the multipole 
  decomposition to enable further analysis of this representation.
Until now, the meaning of the individual deviators is, to the best of our
  knowledge, not known for cases higher than order two.
This paper should enable the analysis of higher-order tensors
  using the deviatoric decomposition within a single framework.
To gain a first impression we analyze the second-order multipole
  decomposition for symmetric tensors in more detail and discuss a 
  connection between the multipoles and the eigenvectors.
  
The paper is organized in the following way.
First, a general introduction to the required tensor algebra is given.
Then, the mathematical background of the deviatoric decomposition 
  and the multipole representation is presented.
Because of their prevalence and familiarity, the second- and  fourth-order decompositions are described in 
  more detail to explain the calculations.
The final chapter describes a more or less well known application 
  of the multipole decomposition:
The multipoles of a stiffness tensor can be used to calculate
  the anisotropy type of the underlying material.
However, as stated earlier, the exact meaning of the various deviators is, to the best
  of our knowledge, not known.
In further works we want to change this.
Identifying the meaning of such abstract object requires interpretations from a range of different fields which can all contribute unique insights. 
Therefore, we want to provide other
  researchers with the necessary tools and concepts to analyze the meaning of the deviators.
Furthermore, we hope to find more applications for the deviatoric 
  decomposition and the multipole representation to learn more 
  about their potential meaning in the light of their generality.
\section{Tensor algebra}
% Euclidean space
The reader may consider the dimension $n$ to be three in most cases.
To ease the description without loss of generality, we will assume an orthonormal
  basis at all times, so there will be no distinction between
  co- and contravariant tensors and all indices will be lower indices.

%Euclidean space, inner product, orthonormal basis, Kronecker delta
We denote the $n$-dimensional \textbf{Euclidean vector space} 
  as $\mathcal{V}_n$.
Its \textbf{scalar product} is a bilinear mapping of two vectors 
  $\textbf{x}, \textbf{y}$ to a real number and is denoted as 
  $\textbf{x} \cdot \textbf{y}$.
We assume that we are given an orthonormal basis 
  $\{\textbf{e}_1,\ldots, \textbf{e}_n\}$, i.e.
  $\textbf{e}_i \cdot \textbf{e}_j = \delta_{ij}$
with the Kronecker delta $\delta_{ij}$.
% tensor of order o as multilinear map
Because no distinction between co- and contravariant tensors is required, 
we define an $n$-dimensional \textbf{tensor} of order $q$ as 
  multilinear map of $q$ vectors to the real numbers
\begin{equation}
    \mathds{T} : (\mathcal{V}_n)^q \to \mathds{R} \text{.}
\end{equation}
As a multilinear map, the tensor can also be described by its 
  coefficients with respect to a fixed orthonormal basis of 
  $\mathcal{V}_n$, say $\{\textbf{e}_i\}$
\begin{equation}
    \mathds{T}(\textbf{e}_{i_1},\ldots,\textbf{e}_{i_q}) = T_{i_1,\ldots,i_q}
\end{equation}
Therefore, a zeroth-order tensor can be represented as a scalar, 
  a first-order tensor as a vector of dimension $n$ and 
  a second-order tensor as an $n \times n$-matrix. 
Higher-order tensors can be represented as arrays of order $q$. 
For example, the coefficients of a three-dimensional fourth-order tensor 
  can be represented as $3\times 3\times 3\times 3$ array. 
In some cases, the tensors in this paper will be described by this index
  notation.
Then the number of the indices describes the order of the tensor.
If there is another orthonormal basis 
  $\{\textbf{f}_1, \ldots, \textbf{f}_n\}$
  of $\mathcal{V}_n$ expressed in the basis $\{\textbf{e}_i\}$ as
  $\textbf{f}_j = \sum_{i=1}^n c_{ij} \textbf{e}_i$
  then $\mathds{T}$ is represented by the numbers
\begin{equation}
  T'_{j_1 j_2 \ldots j_q} = \sum_{i_1=1}^n \ldots \sum_{i_q=1}^n c_{i_1 j_1} \ldots c_{i_q j_q} T_{i_1 \ldots i_q} 
\end{equation}
  with respect to the basis $\{\textbf{f}_j\}$.
In the following, we describe our tensors always as coefficient arrays of
  order $q$, i.e. as represented with respect to the basis $\{\mathbf{e}_i\}$.
This is also the representation in our (and probably any typical) 
  visualization application.

% addition and scalar multiplication create vector space
Let $\mathds{S}, \mathds{T}$ be tensors of order $q$ over 
  $\mathcal{V}_n$, and $a,b\in \mathds{R}$. 
Then, the set of all $q$-arrays, with the scalar multiplication 
  and the tensor addition
\begin{equation}
    a\mathds{S} + b \mathds{T} = aS_{i_1 i_2 \dots i_q} + bT_{i_1 i_2 
    \dots i_q}
\end{equation}
  create a $n^q$-dimensional vector space $\mathcal{V}_n^q$ representing 
  all $n$-dimensional tensors of order $q$.

% tensor product
An often used tensor operation is the \textbf{tensor product} or
  \textbf{outer product}. 
  The tensor product of a $q$th-order tensor $\mathds{A}$ and 
  a $r$th-order tensor $\mathds{B}$ results in a $(q+r)$th-order tensor 
  $\mathds{C}$ as follows
\begin{equation}
  \mathds{C}_{i_1 \dots i_q j_1 \dots j_r} = \mathds{A}\otimes \mathds{B} = A_{i_1 \dots i_q} B_{j_1 \dots j_r} 
\end{equation}
Mathematically, this creates an algebra of all tensors of arbitrary order 
  over $\mathcal{V}_n$.
\begin{comment}
Another arithmetic operation will be marked by '$\odot$' and is given by
\begin{equation}
    \mathds{A} \odot \mathds{B} = \frac{1}{2} \left[ (\mathds{A} \otimes \mathds{B}^T)^{\stackrel{23}{T}} + (\mathds{A} \otimes \mathds{B})^{\stackrel{24}{T}} \right]
\end{equation}
where {\tiny$\stackrel{ij}{T}$} is the transposition of the $i$th and the $j$th base vector.
\end{comment}

% contraction
Another tensor operation we use is called \textbf{tensor contraction}. 
The tensor contraction is the summation over a determined number 
  of indices. 
The \textbf{single contraction} 
  $ \mathds{C} = \mathds{A} \cdot \mathds{B}$ is 
  the summation of two tensors $\mathds{A}$ and $\mathds{B}$ over
  one index
\begin{equation}
  \begin{split}
     C_{i_1 i_2 \dots i_{n-1} j_2 \dots j_m } &= \mathds{A} \cdot \mathds{B} = \sum\limits_{k=1}^{n} A_{i_1 i_2 \dots i_{n-1} k} B_{k j_2 \dots j_m} \\
     &= A_{i_1 i_2 \dots i_{n-1} k} B_{k j_2 \dots j_m} 
  \end{split}
\end{equation}
where the implicit summation over a repeated index, i.e. $A_{i_1 i_2 \dots i_{n-1} k} B_{k j_2 \dots j_m} $, is called Einstein summation convention. 
The \textbf{double contraction} $\mathds{A}: \mathds{B}$ is 
  analogously defined as 
\begin{equation}
\begin{split}
C_{i_1 i_2 \dots i_{n-2} j_3 \dots j_m } &= \mathds{A} : \mathds{B} = \sum\limits_{k=1}^{n} \sum\limits_{l=1}^{n} A_{i_1 i_2 \dots i_{n-2} k l} B_{kl j_3 \dots j_m} \\
&= A_{i_1 i_2 \dots i_{n-2} kl} B_{kl j_3 \dots j_m} 
\end{split}
\end{equation}

% trace, determinant
One can define a trace $\text{tr}_{i,j}$ for each index pair $i,j$,
  but here we will define the general \textbf{trace} $\text{tr}()$ of a tensor
  as the following sum about the first two indices
\begin{equation}
    \text{tr } ( \textbf{T} ) = \text{tr}_{1,2} ( \textbf{T} ) = \sum\limits_{s=1}^n D_{ss i_3 i_4 \dots i_n}
\end{equation}
  which generalizes the notion well-known for second-order tensors.
Therefore, a tensor of order higher than two has more than one trace.
The \textbf{determinant} of a second-order tensor $\textbf{T}$ is 
  defined as
\begin{align}
\det (\textbf{T}) \quad=\quad & T_{11}T_{22}T_{33} + T_{12}T_{23}T_{13} + T_{13}T_{12}T_{23} \\
& -T_{13}^2 T_{22} - T_{23}^2 T_{11} -T_{12}^2 T_{33} \nonumber
\end{align}

Let $\mathcal{S}_q$ be the $q$\textsuperscript{th}-order symmetric group of all 
  $q!$ permutations $\pi$ on the integers $(1, 2, \dots, q)$.
Let $\mathds{T}$ be a three-dimensional tensor of order $q$. 
Then the permutation of the tensor $\pi \mathds{T}$ can be 
  defined as 
\begin{equation}
    (\pi \mathds{T})_{i_{\pi_1} \dots i_{\pi_q}} = T_{i_1 
    \dots i_q}
\end{equation}
% symmetries (minor symmetries, major symmetry, total symmetry)
We call a tensor $\mathds{T}$ \textbf{totally symmetric}
  if $\pi \mathds{T} = \mathds{T}$ for all 
  $\pi \in \mathcal{S}_q$. 
So, every symmetric second-order tensor is totally symmetric.
A tensor of order $2n$ can have the so called \textbf{minor symmetry},
  which means it is invariant to index permutations in the sets 
  $\{1,2\}$ and $\{2n-1,2n\}$ of index positions.
These even-order tensors can also have the \textbf{major symmetry}, 
  which describes the invariance of a permutation of the index positions
  $\{ 1,\dots , n \}$ with $\{ n+1, \dots ,2n\}$.
For a fourth-order three-dimensional tensor $\mathds{T}$
  the minor symmetries and the major symmetry are given by
\begin{equation}
    T_{ijkl} = T_{jikl} = T_{ijlk}, \qquad T_{ijkl} = T_{klij}.
\end{equation}
It is important for the following discussion that total symmetry
  implies more symmetry than minor and major symmetries 
  alone.
This becomes clear if one realizes that a totally symmetric
  three-dimensional fourth-order tensor has $15$ independent
  coefficients, i.e. $\dim(\mathcal{S}^4)=15$.
In contrast, a three-dimensional fourth-order tensor with
  minor and major symmetries has $21$ independent 
  coefficients.
Based on these definitions, we define the \textbf{totally 
  symmetric part} $\textbf{s}\mathds{T}$ of a general tensor as
\begin{equation}
    \textbf{s}\mathds{T} = \frac{1}{q!} \sum\limits_{\pi \in 
    \mathcal{S}_q} \pi \mathds{T}
\end{equation}
% antisymmetric part
We call the remainder 
\begin{equation}
    \textbf{a}\mathds{T} = \mathds{T} - \textbf{s}\mathds{T}
\end{equation}
  the \textbf{asymmetric part} of $\mathds{T}$.
Obviously, the asymmetric part of any tensor has no 
  totally symmetric part which creates a vector space
  of all asymmetric tensors which is actually orthogonal to
  the totally symmetric tensor space, see~\cite{backus}.
The symmetric and traceless part of a tensor $\mathds{T}$ will be
  described by $\lfloor \mathds{T} \rfloor $.
For the tensor product $\textbf{n}\otimes \textbf{m}$ it is
\begin{equation*}
    \lfloor \textbf{n} \otimes \textbf{m} \rfloor = \frac{1}{2} \left( \textbf{n} \otimes \textbf{m} + \textbf{m} \otimes \textbf{n} \right) - \frac{1}{3} \left(\textbf{n} \cdot \textbf{m}\right) \textbf{I}.
\end{equation*}

% deviator (=totally symmetric and traceless) (mention the coincidence with second order case
This work focuses on so called \textbf{deviators}.
In the visualization literature and classical mechanics 
texts, a deviator is defined as a traceless tensor of   
  second-order.
Here, we extend this definition following Zou et 
  al.~\cite{zou2001} to general three-dimensional 
  $q$th-order deviators. 
A \textbf{deviator} $\mathds{D}^{(q)}$ is a 
  three-dimensional tensor 
  of arbitrary order $q$ that is totally symmetric and 
  traceless, i.e.
\begin{equation}
  D_{i_1 i_2 i_3 \dots i_q} = D_{i_2 i_1 i_3 \dots i_q} 
    = \dots = D_{i_q i_2 i_3 \dots i_1}, 
    \qquad \text{tr }\left( \mathds{D}^{(q)}\right) = 0 
\end{equation}
\section{Spectral Decomposition}
The best known decomposition for $n\times n$ tensors is the 
  spectral decomposition.
%Second-order spectral decomposition
Such a tensor $\textbf{T}$ can be described by three $n$-dimensional
  eigenvectors $\textbf{v}_i$ and three eigenvalues $\lambda_i$
\begin{equation}
    \textbf{T} = \sum\limits_{i=1}^3 \lambda_i \textbf{v}_i \otimes \textbf{v}_i
\end{equation}
For symmetric tensors the results are real.
  
%Fourth-order spectral decomposition
General tensors of order higher than two can not be decomposed
  using this decomposition.
There are some special cases where a generalization of the spectral
  decomposition exists.
One of these special cases is for example a fourth-order tensor
  $\mathcal{T}$ with the two minor symmetries
  and the major symmetry.
By choosing a suitably defined tensor basis exploiting these symmetries, this tensor can be mapped onto a symmetric second-order six-dimensional 
  tensor using the \textbf{Mandel} or \textbf{Kelvin mapping}  \cite{Kelvin:1856,Mehrabadi1990,Nagel2016} as follows:
%\begin{equation}
%    11\to 1, \quad 22\to 2, \quad 33\to 3, \quad 23 \to 4, \quad 13 \to 5, \quad 12 \to 6
%\end{equation}
\begin{equation}
    K_{mn}(T_{ijkl}) = \begin{pmatrix}
    T_{1111} & T_{1122} & T_{1133} & \sqrt{2}T_{1123} & \sqrt{2}T_{1113} & \sqrt{2}T_{1112}\\
    T_{1122} & T_{2222} & T_{2233} & \sqrt{2}T_{2223} & \sqrt{2}T_{2213} & \sqrt{2}T_{2212}\\
    T_{1133} & T_{2233} & T_{3333} & \sqrt{2}T_{3323} & \sqrt{2}T_{3313} & \sqrt{2}T_{3312}\\
    \sqrt{2}T_{1123} & \sqrt{2}T_{2223} & \sqrt{2}T_{3323} & 2T_{2323} & 2T_{2313} & 2T_{2312}\\
    \sqrt{2}T_{1113} & \sqrt{2}T_{2213} & \sqrt{2}T_{3313} & 2T_{2313} & 2T_{1313} & 2T_{1312}\\
    \sqrt{2}T_{1112} & \sqrt{2}T_{2212} & \sqrt{2}T_{3312} & 2T_{2312} & 2T_{1312} & 2T_{1212}
    \end{pmatrix}.
\end{equation}
This tensor can be decomposed by the known second-order spectral
  decomposition with the eigenvalues $\Lambda_i$.
The results are six first-order six-dimensional tensors.
These can be mapped using the inverse Kelvin mapping onto second-order
  three-dimensional tensors $M_i$.
These are the eigentensors of the fourth-order tensor $\mathcal{T}$, which
  can be represented by
\begin{equation}
    \mathcal{T} = \sum\limits_{i=1}^6 \Lambda_i  \mathbf{M}_i\otimes \mathbf{M}_i.
\end{equation}

\section{Spherical Harmonics}
There is a well known relation between totally symmetric tensors
  and spherical harmonics.
This section will describe this connection, following
the representation in Backus' \cite{backus} paper. 
% description follows seminal paper by Backus (1970)
% permutation of a tensor

\begin{comment}
%Brauche ich glaube ich nicht
Then the $n$th-order deviator is given by
\begin{equation} \label{dev}
\textbf{D}^{(n)} = c_1^{(n)} \textbf{P}^{(n)} + c_2^{(n)} \textbf{Q}^{(n)}, \qquad c^{(n)}= c_1^{(n)} + i c_2^{(n)}
\end{equation}
with $c_1^{(n)}, \, c_2^{(n)} \in \mathds{R}$ and $\textbf{P}^{(n)}$ and $\textbf{Q}^{(n)}$ being the real and imaginary parts of $\textbf{w}^{\otimes n} = \underbrace{\textbf{w}\otimes \textbf{w} \otimes \dots \otimes \textbf{w}}_{n \text{times}}$ where $\textbf{w} = \textbf{e}_1 + i\textbf{e}_2$ with the orthogonal frame $\{ e_i\}$ of the second-order euclidean space.
\end{comment}

% homogeneous polynomials are isomorphic to totally symmetric tensors
Let $\mathds{T}$ be a three-dimensional tensor of order $q$.
We consider the polynomial
\begin{equation}\label{polynom}
    \mathds{T}(\textbf{r}) = T_{i_1 \dots i_q} r_{i_1} \dots 
    r_{i_q}
\end{equation}
  where $\textbf{r}=(r_1,r_2,r_3)$ are coordinates with 
  respect to our orthonormal basis 
  $(\mathbf{e}_1,\mathbf{e}_2,\mathbf{e}_3)$.
All monomials of this polynomial have the degree 
  $q$, so $\mathds{T}(\textbf{r})$ is a 
  homogeneous polynomial of degree $q$.
It is called the \textbf{polynomial generated by} 
  $\mathds{T}$. 
Since the products of the $r_{i_j}$ are products of real 
  numbers and therefore commutative, we have 
  $\mathds{T}(\textbf{r}) = \mathds{S}(\textbf{r})$ 
  if $\mathds{S}$ is the totally symmetric part of 
  $\mathds{T}$. 
More specifically, it is
\begin{equation}
    S_{i_1 \dots i_q} = \frac{1}{q!} \partial_{i_1} \dots 
    \partial_{i_q} [T(\textbf{r})]
\end{equation}
Therefore, two different totally symmetric tensors generate 
  different polynomials, but two tensors with the same
  totally symmetric part generate the same polynomial. 
Let $\mathcal{P}^q$ be the linear space of all homogeneous 
  polynomials of degree $q$ in three dimensions. 
Then, we have an isomorphism between $\mathcal{P}^q$ and $\mathcal{S}^q$,
  the space of the totally symmetric tensors of order $q$. 

% harmonic polynomials equal the traceless totally symmetric tensors
Harmonic polynomials $\mathds{P}^{(q)}(\mathbf{r})$ in 
  $\mathcal{P}^q$ are called \textbf{spherical harmonics}. 
The well-known decomposition of homogeneous polynomials, 
  interpreted as functions on the sphere, by spherical harmonics 
  states the following:
For every polynomial in $\mathcal{P}^q$,
  there are spherical harmonics $\mathds{H}^{(q)} (\textbf{r}), 
  \mathds{H}^{(q-2)}(\textbf{r}), \dots $ such that
\begin{equation}\label{spherical}
    \mathds{P}^{(q)} (\textbf{r}) = \mathds{H}^{(q)}(\textbf{r}) + \textbf{r}^2 \mathds{H}^{(q-2)} 
    (\textbf{r}) + \dots
\end{equation}
Applying the three-dimensional Laplacian, it follows that 
  $\mathds{H}^{(q)}$ is unique. 
%\eqref{spherical} can be represented as
%\begin{equation}
%    P^{(q)}(\textbf{r}) = H^{(q)}(\textbf{r}) + 
%       \textbf{r}^2 P^{(q-2)}(\textbf{r})
%\end{equation}
%with a unique $P^{(q-2)} \in \mathcal{P}^{q-2}$.
This allows a nice characterization of deviators.
Let $\mathds{H}$ be a $q$\textsuperscript{th}-order totally symmetric tensor 
  that generates a harmonical polynomial 
  $\mathds{H}(\textbf{r})$. 
It follows $\Delta \mathds{H}(\textbf{r})=0$ and therefore
\begin{equation}
    H_{jji_3 \dots i_q} r_{i_3} \dots r_{i_q} = 0
\end{equation}
Therefore, such a tensor is traceless and totally symmetric,
  i.e. a deviator.
% decomposition of spherical harmonics allows to decompose totally symmetric tensors
For a general totally symmetric tensor $\mathds{S}$ of 
  order $q$, we may use the spherical harmonic decomposition
  of $\mathds{S}(\mathbf{r})$ to compute tensors $\mathds{H}^{(q)}$ 
  using \eqref{spherical}. 
By comparing the coefficients, we can derive that
\begin{equation}\label{deviatoric}
    \mathds{S} = \mathds{H}^{(q)} + \textbf{s} 
    (\textbf{I}\mathds{H}^{(q-2)}) + \dots
\end{equation}
  where $\mathbf{I}$ is the identity tensor and $\textbf{s}()$
  denotes the totally symmetric part of a tensor as before.
Therefore, any totally symmetric tensor in three dimensions
  is uniquely described by a series of deviators $\mathds{H}_i$, which will 
  be described in the following by $\mathds{D}_i$, of order 
  $q, q-2, \ldots$.

% deviators as irreducible components
% antisymmetric tensors in 3D allow for isomorphic mapping to 
% symmetric tensors two orders lower

% multipole description of deviators according to Hamilton 
% and Sylvester

\section{Deviatoric Decomposition}
As mentioned before, the connection between totally symmetric 
  tensors and spherical harmonics is well known.
Thus, each totally symmetric tensor can be described by
  deviators.
Also, every tensor of arbitrary order up to dimension three can be 
  described by deviators uniquely.
Backus \cite{backus} described this so-called \textbf{deviatoric
  decomposition}.

It is well known that every tensor $\mathds{T}$ of order $q$ can be
  decomposed into a totally symmetric $\mathbf{s}\mathds{T}$ and an asymmetric
  part
  $\mathbf{a}\mathds{T}$
  \begin{equation*}
      \mathds{T} = \textbf{s}\mathds{T} + \textbf{a}\mathds{T} 
  \end{equation*}

The deviatoric decomposition of the totally symmetric part is 
  described in the previous section using the spherical harmonics.
The deviatoric decomposition of the asymmetric part is
  far less well-known.
The gap is bridged by a unique isomorphism from the totally symmetric tensor 
  space $\mathcal{S}^p$ into the asymmetric tensor space 
  $\mathcal{A}^q$ of order $q$
  \begin{equation*}
      \phi ( \mathcal{S}^p ) \to \mathcal{A}^q.
  \end{equation*}
Through Schur's Lemma \cite{wigner} it is known that the space of the isormorphisms uniquely determines
  the isomorphism.
These totally symmetric tensors can then be decomposed into 
  deviators.
Thus, eventually this isomorphism allows the decomposition of the asymmetric part of the
  tensor into deviators.

\begin{figure}[h]
    \centering
    \includegraphics[scale = 0.22]{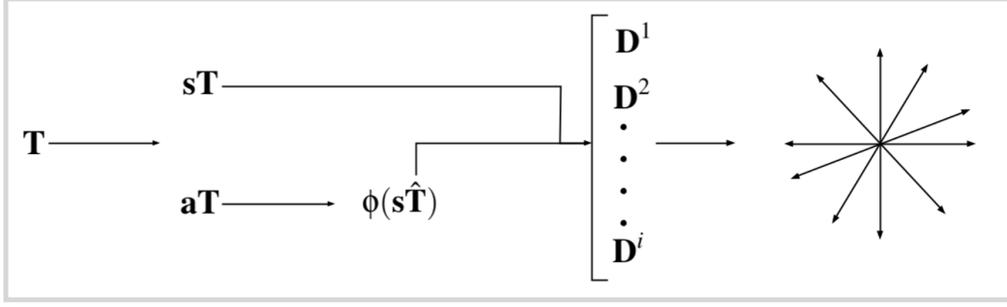}
    \caption{Each tensor $\textbf{T}$ of every order up to dimension three can be decomposed into a set of vectors (of different length or a set of vectors and scalars). First, it must be decomposed into its symmetric part $\textbf{sT}$ and its asymmetric part $\textbf{aT}$. For the asymmetric part an isomorphism can be used to represent this part also with a symmetric tensor $\textbf{s}\hat{\textbf{T}}$. These symmetric tensors can be decomposed in deviators $\textbf{D}^j$. Each of these deviators can be witten as the symmetrization of a set of vectors, called multipoles, multiplied by a scalar.}
    \label{fig:aniso}
\end{figure}

\section{Multipoles}
We complete our description by decomposing deviators
  even further as tensor products of vectors which
  are called \textbf{multipoles}.
This result was first described by 
  Maxwell~\cite{Maxwell:1881}, who found it as an elegant
  geometric description of spherical harmonics.
Using a theorem by Sylvester~\cite{Sylvester:1876}, Maxwell
  stated in our terms that for any three-dimensional deviator
  $\mathds{D}$ of order $q$, there are $q$ unit real vectors
  $\textbf{n}_1, \dots , \textbf{n}_q$ called multipoles
  and a real non-negative number $a$ such that
\begin{equation} \label{eq:multipole}
    \mathds{D} = a \, \lfloor\textbf{n}_1 \otimes \dots \otimes 
    \textbf{n}_q \rfloor
\end{equation}
The unit vectors $\textbf{n}_i$ are uniquely defined by 
  $\mathds{D}$ up to an even number of sign changes, and 
  the number $a$ is unique.

The major algorithmic task is to actually compute the multipoles 
  of all deviators in the deviatoric decomposition of a given
  tensor.
To identify multipoles of a $q$\textsuperscript{th}-order deviator, 
  Zou et al.~\cite{zou} define a polynomial of degree $2q$ over the 
  complex numbers by
\begin{equation} \label{eq:poly}
    a_{q,0} + \sum\limits_{r=1}^{q} \sqrt{\frac{q!q!}{(q+r)!(q-r)!}} 
    \left( x^r \overline{a_{q,r}} + (-1)^r x^{-r} a_{q,r}\right) = 0 .
\end{equation}
We do not explain the background of this polynomial here, but refer the reader interested in such detail to
  the original paper by Zou et al~\cite{Zou:2003}.
There, one can see that the polynomial is based on the construction 
  of an orthogonal basis of the corresponding space of deviators 
  for each $x \in \mathds{C}$ and each order $q$.
The coefficients $a_{q,r}$, $r=0,\ldots,q$ of the polynomial
  are the coefficients of the deviator expressed in terms of the basis.
The idea is that a root of the polynomial selects a basis that allows
  one to read off the multipoles.
The coefficients $a_{q,r}$ for an explicit example are given in
  \autoref{sec:MultipoleDecomposition}.
However, the roots of the above polynomial come in pairs as
  $\overline{x}_r^{-1}$ is a root if $x_r$ is a root.
Therefore, we actually get $2\cdot q$ different roots, but we need
  only one member of each pair to compute the multipoles.
The multipoles of our deviators are found via
\begin{equation} \label{eq:zeros}
    x_r = e^{i\varphi_r} \tan \frac{\theta_r}{2}
\end{equation}
\begin{equation}\label{eq:mps}
    \textbf{n}_r = 
    \textbf{e}_3 \cos(\theta_r) + (\textbf{e}_1 \cos 
    \varphi_r + \textbf{e}_2 \sin \varphi_r) \sin \theta_r
    \quad
    r=1,\dots 8.
\end{equation}

\section{Second-order tensor decomposition}
Based on the general considerations above, we now study the 
important special case of a second-order tensor to apply the method
on familiar ground.
A general second-order tensor can be decomposed into 
  the zeroth-order deviator $d$, the first-order 
  deviator $\textbf{d}$ and the second-order deviator
  $\textbf{D}$ by
\begin{equation}
         T_{ij} = d I_{ij} + \varepsilon_{ijk} \cdot d_{k} + D_{ij}, \quad d = \frac{1}{3} T_{ii}, \quad d_i = \frac{1}{2} \varepsilon_{ijk} T_{jk}
     \end{equation}
%
\begin{comment}  
\notegs{Irgendwie gefällt mir unsere Notation nicht. 
        Ich fände es besser, wenn die einzelnen Deviatoren
        den gleichen Buchstaben wie der ursprüngliche Tensor
        verwenden, zusammen mit diversen Ergänzungen
        (hochgestellte Zahlen, Dach, Schlange etc. um die
        Ordnung und die Nummer anzugeben, wenn es mehrere
        Deviatoren gleicher Ordnung in der Zerlegung gibt.
        Andernfalls ist man total verwirrt, wenn mehrere
        Tensoren zerlegt werden und zudem womöglich viele
        verschiedene tensorielle Größen in den physikalischen
        Gleichungen verwendet werden.}
\end{comment}
%
$\bm{\varepsilon}$ describes the third-order permutation
  tensor and $\textbf{I}$ the second-order identity tensor.
If $\mathbf{T}$ is symmetric the antisymmetric part represented by
  $\textbf{d}$ becomes zero.
  
\subsection{Relation to the Eigendecomposition} \label{sec:eigen}
In the second-order three-dimensional symmetric tensor case,
  the tensor $\textbf{T}$ can be decomposed by the
  spectral decomposition into three vectors $\textbf{v}_i$ and three scalars $\lambda_i$
  \begin{equation}
      \textbf{T}=\sum\limits_{i=1}^3\lambda_i \textbf{v}_i \otimes \textbf{v}_i
  \end{equation}
  or by the deviatoric decomposition into two vectors $\textbf{m}_1$
  and $\textbf{m}_2$ and two scalars $d$ and $a$
  \begin{equation}
      \textbf{T}= d\cdot \textbf{I} + a \lfloor \textbf{m}_1 \otimes \textbf{m}_2 \rfloor
  \end{equation}
This suggests there is a connection between the two decompositions.
This connection is analyzed in the sequel.
\begin{comment}
\notegs{Wäre es nicht sinnvoll, die beiden Zerlegungen als
  Formeln direkt gegenüberzustellen? Also $T=\lambda_i v_i \otimes v_i$ für die Spektralzerlegung und die andere Formel für die 
  deviatorische Zerlegung?}
\end{comment}
On the one hand, there are three orthonormal eigenvectors and 
  three real eigenvalues that describe such a tensor.
On the other hand, there are two multipoles $\textbf{m}_1$ and 
  $\textbf{m}_2$ that describe the second-order deviator of
  this tensor.
Together with the isotropic part, i.e. trace times identity, this
  is the deviatoric decomposition of this tensor.
However, because eigendecomposition, deviatoric decomposition and multipole
  decomposition of deviators are all unique,
  there must be a connection.
This will be elaborated on in this section.

There are three different cases for the position of two multipoles.

\textbf{Case 1.} The scalar $a$ from Eq.~\eqref{eq:multipole} equals zero
  or the multipoles are given by
  $\textbf{m}_1 = \textbf{m}_2 = \left( 0, 0, 0 \right)^\text{T}$.
This is the case, if and only if the tensor has a triple eigenvalue, i.e. is an isotropic or spherical tensor.

\textbf{Case 2.} The multipoles $\textbf{m}_1$ and $\textbf{m}_2$
  are identical, if and only if the tensor has a double eigenvalue.
Then, the multipoles are given by the eigenvector according to the 
  eigenvalue with the largest absolute value.
 % \notegs{Stimmt das mit dem betragsmäßig größtem Eigenwert? Ist es nicht der Eigenvektor zum einfachen Eigenwert? - Chiara: Siehe Beweis Idee.}

\texttt{Proof idea.}
Assume $\textbf{m}_1=\textbf{m}_2 =: \textbf{m}$. Then 
\begin{equation}
  \textbf{T}=\text{tr}(\textbf{T})\textbf{I} + a \lfloor \textbf{m} \otimes \textbf{m} \rfloor
\end{equation}
We know that adding a multiple of the identity to a tensor does 
  not change the
  eigenvectors of a tensor and does not change the multiplicity of the eigenvalues,
  so we define $\textbf{m}_1 \cdot \textbf{m}_2 =:b$.
Because the multipoles are the same, the symmetrization of the tensor product is
  the tensor product itself and reads
\begin{equation}
    a\lfloor \textbf{n} \otimes \textbf{n} \rfloor = \underbrace{a\begin{pmatrix}
    m_1^2 & m_1 m_2 & m_1 m_3\\
    m_1 m_2 & m_2^2 & m_2 m_3\\
    m_1 m_3 & m_2 m_3 & m_3^2
    \end{pmatrix}}_{=: \textbf{M}} - ab \textbf{I}
\end{equation}
The eigenvalues of $\textbf{M}$ are given by
\begin{equation}
    \lambda_1 = 0, \quad \lambda_2 = 0, \quad \lambda_3 = a \cdot (m_1^2 + m_2^2 + m_3^2).
\end{equation}
Thus, the tensor has a double eigenvalue.

Assume now, the tensor is degenerated and $s$ is the double
  eigenvalue.
Zheng \cite{zheng2006degenerate} states that a degenerate 
  second-order three-dimensional tensor can be represented by
\begin{equation}
    T = s\textbf{I} \pm \textbf{V}\cdot \textbf{V}^\text{T}
\end{equation}
  where $\textbf{V}$ is defined by the eigenvector to the single eigenvalue by
  $\textbf{e} = \textbf{V}/\left\lVert\textbf{V}\right\rVert$.
It follows
\begin{equation}
    \textbf{V}\textbf{V}^\text{T} = \textbf{e} \left\lVert\textbf{V}\right\rVert \cdot \left(\textbf{e} \left\lVert\textbf{V}\right\rVert\right)^\text{T} = \left\lVert\textbf{V}\right\rVert^2( \textbf{e} \cdot \textbf{e}^\text{T} )
    = \left\lVert\textbf{V}\right\rVert^2\begin{pmatrix}
    e_1^2 & e_1e_2 & e_1e_3\\
    e_1e_2 & e_2^2 & e_2e_3\\
    e_1e_3 & e_2e_3 & e_3e_3
    \end{pmatrix}
\end{equation}
We can derive $\textbf{m}_1 = \textbf{m}_2 = \textbf{e}$. 
Therefore, the multipoles collapse into one vector and equal the eigenvector
  corresponding to the eigenvalue with the largest absolute value, which is
  in this case the single eigenvalue.
Thus, the assumption is true.

\textbf{Case 3.} When the multipoles $\textbf{m}_1$ and $\textbf{m}_2$ 
  are different,
 a connection between the eigenvectors and eigenvalues 
  of the symmetric part and the multipoles can be found.
The eigenvector deduced by the, according to the absolute 
  value, largest eigenvalue and the, sorted by absolute value, 
  medium eigenvalue are the bisecting lines of the two multipoles.
The eigenvector associated to the eigenvalue with largest 
  absolute value is the bisecting line which has the smaller 
  angle to the multipoles.
  In other words, the closer two of the eigenvalues become, the closer the
  multipoles bias towards the direction of the largest (absolute) eigenvalue.
In the eigenvector system where $\left| \lambda_1 \right| > 
  \left| \lambda_2 
  \right| > \left| \lambda_3 \right|$, the multipoles are given by 
\begin{equation}\label{eq:multiDifferentEvecs}
    \mathbf{m}_1 = \pm\begin{pmatrix}
    \sqrt{\frac{2\lambda_1 +\lambda_2}{a}}
 \\
    \sqrt{\frac{-\lambda_1 -2\lambda_2}{a}}\\
    0   \end{pmatrix}
    \qquad \mathbf{m}_2 = \pm\begin{pmatrix}
    \sqrt{\frac{2\lambda_1 + \lambda_2}{a}}\\
    -\sqrt{\frac{-\lambda_1-2\lambda_2}{a}}\\
    0
    \end{pmatrix}
\end{equation}
The angle between the first multipole and the eigenvector 
  $\textbf{e}_1$ is given by
\begin{equation}
    \arccos \left( \frac{\sqrt{\frac{2\lambda_1 +\lambda_2}{a}}}    
                        {\sqrt{\frac{\lambda_1 - \lambda_2}{a}}}\right).
\end{equation}

\texttt{Proof idea.} 
Let $\left| \lambda_1 \right| > \left| \lambda_2 \right| > 
  \left| \lambda_3 \right|$, where $\lambda_i$ are the eigenvalues of the 
  deviator.
Because the deviator is traceless, $\lambda_3 = -(\lambda_1 + \lambda_2)$ holds.
The idea for proving the other cases is analogous.
Assume, the eigenvectors are given by the coordinate axes 
  $\textbf{e}_1, \, \textbf{e}_2, \, \textbf{e}_3$.
Then the deviator $\textbf{D}$ is given by
\begin{equation*}
    \textbf{D} = \text{diag}\, ( \lambda_1, \lambda_2, \lambda_3 )
\end{equation*}
This tensor can also be written as
{\small
\begin{equation}
\begin{split}
    \textbf{D} &= a \begin{pmatrix}
    \sqrt{\frac{2\lambda_1+ \lambda_2}{a}}^2 & 0 & 0\\
    0 & -\sqrt{\frac{-\lambda_1 + 2\lambda_2}{a}}^2 & 0 \\
    0 & 0 & 0
    \end{pmatrix}
    - \frac{a}{3} \left( \frac{2\lambda_1 + \lambda_2}{a} - \frac{-\lambda_1-2\lambda_2}{a}\right)\textbf{I}
    \end{split}
\end{equation}}
Using \eqref{eq:multiDifferentEvecs} it follows
\begin{equation*}
    \textbf{D} = a \cdot \lfloor \textbf{m}_1 \otimes \textbf{m}_2 \rfloor
\end{equation*}
Then, the rotation invariance must be proven and the result follows through the equivalence of 
  each step.

This connection underlines the importance of the multipoles.
We do not claim that the multipole decomposition is better than the
  eigendecomposition.
We just want to demonstrate that in the symmetric second-order 
  case, the information in the eigendecomposition can also 
  be found in the multipole decomposition.
Thus, an analysis of higher-order tensors by using the deviatoric
  decomposition appears reasonable.
With this connection, we want to emphasize the close connection 
  between the two decompositions to show the importance of 
  multipoles, in particular as they allow a generalization
  to tensors not amenable to an eigendecomposition.
  
\section{Fourth-order tensor decomposition}
The deviatoric decomposition described above can be used for every
  tensor of any order up to dimension three.
Backus \cite{backus} gave a decomposition of a general fourth-order
  three-dimensional tensor.
It can be decomposed into a fourth-order deviator 
  $\mathcal{D}$, three third-order deviators $\mathrm{D}^{(i)}$, six 
  second-order deviators $D^{(i)}$, six first-order deviators $\mathbf{d}^{(i)}$
  and three zeroth-order ones $d^{(i)}$.
The tensor is then given by 
\begin{equation}
\begin{split}
    \mathcal{T}_{ijkl} = & D_{ijkl} \\
    &+ \varepsilon_{ikm}D^{(1)}_{jlm} + \varepsilon_{jkm}D^{(1)}_{ilm} + \varepsilon_{ilm} D^{(1)}_{jkm} + \varepsilon_{jlm} D^{(1)}_{ikm}\\
    &+ D^{(2)}_{ijm}\varepsilon_{mkl}\\
    &+ D^{(3)}_{klm}\varepsilon_{mij}\\
    &+ \delta_{ij} D^{(1)}_{kl} + \delta_{kl}D^{(1)}_{ij} + \delta_{ik} D^{(1)}_{jl} + \delta_{jl} D^{(1)}_{ik} + \delta_{il} D^{(1)}_{jk} + \delta_{jk} D^{(1)}_{il} \\
    &+ \delta_{ij} D^{(2)}_{kl} + \delta_{kl} D^{(2)}_{ij} - \frac{1}{2} \delta_{ik}D^{(2)}_{jl} - \frac{1}{2} \delta_{jl} D^{(2)}_{ik} - \frac{1}{2} \delta_{il} D^{(2)}_{jk} - \frac{1}{2} \delta_{jk} D^{(2)}_{il}\\
    &+\delta_{ij}D^{(3)}_{kl}-\delta_{kl}D^{(3)}_{ij}\\
    &+D^{(4)}_{ik}\delta_{jl} + D^{(4)}_{jk}\delta_{il} - D^{(4)}_{il}\delta_{jk} - D^{(4)}_{jl}\delta_{ik}\\
    &+D^{(5)}_{ki}\delta_{lj} + D^{(5)}_{li}\delta_{kj} - D^{(5)}_{kj}\delta_{lj} - D^{(5)}_{lj}\delta_{kj}\\
    &+D^{(6)}_{ik}\delta_{jl} - D^{(6)}_{il}\delta_{jk}\\
    &+(\varepsilon_{ikm}\delta_{jl} + \varepsilon_{jkm}\delta_{il} + \varepsilon_{ilm}\delta_{jk} + \varepsilon_{jlm}\delta_{ik}) d^{(1)}_m\\
    &+d^{(2)}_i \varepsilon_{jkl} + d^{(2)}_j\varepsilon_{ikl} - \frac{2}{3} \delta_{ij}\varepsilon_{klm}d^{(2)}_m\\
    &+d^{(3)}_k \varepsilon_{lij} + d^{(3)}_l\varepsilon_{kij} - \frac{2}{3} \delta_{kl}\varepsilon_{ijm}d^{(3)}_m\\
    &+\delta_{ij}\varepsilon_{klm} d^{(4)}_m\\
    &+\delta_{kl}\varepsilon_{ijm} d^{(5)}_m\\
    &+\varepsilon_{ijk}d^{(6)}_l - \varepsilon_{ijl}d^{(6)}_k - \varepsilon_{kli}d^{(6)}_j + \varepsilon_{klj}d^{(6)}_i\\
    &+d^{(1)} (\delta_{ij}\delta_{kl} + \delta_{ik}\delta_{jl}+\delta_{il}\delta_{jk})\\
    &+ d^{(2)} (\delta_{ij} \delta_{kl} - \frac{1}{2} \delta_{ik}\delta_{jl}-\frac{1}{2} \delta_{il}\delta_{jk} )\\
    &+d^{(3)} (\delta_{ik}\delta_{jl} - \delta_{il}\delta_{jk})
\end{split}
\end{equation}

This decomposition is very complex, but for tensors with some 
  kinds of symmetries, some of the deviators vanish.

\subsection{Stiffness Tensor}
The best known application of the deviatoric decomposition 
  is the calculation of symmetries of materials described by
  the stiffness tensor.
To understand the background of this calculation we give a short
  introduction to the tensor.
  
We describe the stress at a point of the deformed material 
  by the symmetric second-order Cauchy stress tensor 
  $\boldsymbol{\sigma}$ and the strain
  by $\boldsymbol{\varepsilon}$
\begin{equation}
   \boldsymbol{\sigma}= 
     \begin{pmatrix}
      \sigma_{11} & \sigma_{12} & \sigma_{13}\\
      \sigma_{12} & \sigma_{22} & \sigma_{23}\\
      \sigma_{13} & \sigma_{23} & \sigma_{33}
    \end{pmatrix}, \qquad \boldsymbol{\varepsilon} = 
    \begin{pmatrix}
    \varepsilon_{11} & \varepsilon_{12} & \varepsilon_{13}\\
    \varepsilon_{12} & \varepsilon_{22} & \varepsilon_{23}\\
    \varepsilon_{13} & \varepsilon_{23} & \varepsilon_{33}
    \end{pmatrix}
\end{equation} 
The eigenvectors of $\boldsymbol{\sigma}$ are called principal stress directions
  and the eigenvalues principal stresses.
  
\begin{comment}
The gradient of the velocity $\nabla \textbf{v}^\text{T}$ in three
  dimensions is a second-order asymmetric tensor given by 
\begin{equation}
    \textbf{L} = \nabla \textbf{v}^\text{T} \text{ and } [\mathbf{L}]_{ij} = \begin{pmatrix}
    \frac{\partial v_{1}}{\partial v_1} & \frac{\partial v_{2}}{\partial v_1} & \frac{\partial v_{3}}{\partial v_1}\\
    \frac{\partial v_{1}}{\partial v_2} & \frac{\partial v_{v_2}}{\partial 2} & 
    \frac{\partial v_{3}}{\partial v_2}\\
    \frac{\partial v_{1}}{\partial v_3} & \frac{\partial v_{2}}{\partial v_3} & \frac{\partial v_{3}}{\partial v_3}
    \end{pmatrix}
\end{equation}
It describes the spatial change of the velocity field in the current configuration.
\end{comment}

The stiffness tensor describes the linear mapping from strain increments
  into stress increments which can be described by the \textbf{Hooke's law}
\begin{equation}
    \text{d}\boldsymbol{\sigma} = \mathbfcal{C} \mathbf{:} \text{d}\boldsymbol{\varepsilon}.
\end{equation}
It is a fourth-order three-dimensional tensor and can 
  be represented as $3\times 3 \times 3 \times 3$ array.
Assuming a non-polar material such that the Cauchy stress tensor is symmetric,
  and the existence of a 
  scalar potential $\psi$ from which stresses are derived by differentiation
  with respect to a work-conjugate symmetric deformation measure, the stiffness 
  tensor has the two minor symmetries and the major symmetry.
Under these conditions, the number of independent coordinates reduces from
  81 to 21.

\subsubsection{Deviatoric Decomposition}
The deviatoric decomposition of the stiffness tensor can 
  for example be used to compute 
  all symmetry planes for all symmetry classes.
Hergl et al. \cite{hergl2019visualization} used  this fact to visualize the symmetries of the stiffness tensor by designing a glyph.
Since our stiffness tensor is not totally symmetric,
  we also need a decomposition of the asymmetric part
  which was the major achievement of Backus in his paper.
Let $\mathcal{C}$ be a three-dimensional fourth-order tensor
  with major and minor symmetry. 
Then, the totally symmetric and the asymmetric parts are given 
  by
\begin{equation}
    S_{ijkl} = \frac{1}{3} \left( C_{ijkl} + C_{iklj} + 
    C_{iljk}\right), \quad A_{ijkl} = \frac{2}{3} C_{ijkl} - 
    \frac{1}{2} C_{iklj}- \frac{1}{3} C_{iljk}
\end{equation}
The trick, to represent the antisymmetric part by deviators, is to 
  define an isomorphism between the second-order totally symmetric 
  tensors $\mathbf{S}^2$
  and the fourth-order asymmetric tensors 
  $\mathbfcal{A}^4$.
Let $\mathbf{R}$ be a totally symmetric tensor of order two
  in three dimensions.
We define the isomorphism $\phi$ by
\begin{equation}\label{asym}
    \phi(R)_{ijkl} = \delta_{ij} R_{kl} + \delta_{kl}R_{ij} - 
    \frac{1}{2} \delta_{ik}R_{jl} - \frac{1}{2} 
    \delta_{jl}R_{ik} - \frac{1}{2} \delta_{il}R_{jk} 
    -\frac{1}{2} \delta_{jk}R_{il}
    \end{equation}
% this mapping is unique for linear and rotationally invariant
% mappings
This isomorphism is the only possible linear isomorphism 
  that is rotation invariant for any rigid rotation which was proven
  by Backus using the classification of all isotropic tensors 
  by Weyl~\cite{Weyl:1946}. 
% this leads to a decomposition of an arbitrary stiffness 
% tensor into irreducible components
Using \eqref{asym}, the second-order tensor   
  $\mathbf{R}$ is given by  
\begin{equation}
    R_{ij} = \frac{1}{(n-2)} A_{ijll} - 
    \frac{\delta_{ij}}{2(n-1)(n-2)}A_{kkll}
\end{equation}
  and the fourth-order three-dimensional tensor $\mathbfcal{C}$ 
  can be represented by the deviatoric decomposition
\begin{equation} \label{deviatoricDec}
    \mathbfcal{C} = \mathbfcal{D} + 
    6\textbf{s}(\textbf{ID}) + 
    3\textbf{s}(\textbf{II}d) + \phi (\hat{\textbf{D}}) 
    + \frac{1}{2} \phi (\textbf{I}\hat{d}) 
\end{equation}
  where $\textbf{R} = \textbf{D}^2 + \textbf{I} \hat{d} $
  is the deviatoric decomposition of the totally symmetric
  tensor $\textbf{R}$.
Therefore, the stiffness tensor can be uniquely decomposed 
  into one fourth-order deviator $\mathbfcal{D}$, two 
  second-order deviators $\textbf{D}^1$ and $\textbf{D}^2$ 
  and two zeroth-order deviators $d$ and $\hat{d}$.

Finally, let us shortly denote the complete deviatoric 
  decomposition of the stiffness tensor $\mathbfcal{C}$
  in terms of its coefficients which allows to implement it
  even without understanding the theory in this section.
The two zeroth-order deviators are called 
  Lam\'{e} coefficients in engineering and can be computed as
\begin{equation}\label{eq:lame}
  \lambda = d= \frac{1}{15} \left( 2C_{iikk}-C_{ikik} \right),  
  \qquad 
  \mu = \hat{d} = \frac{1}{30} \left( 3C_{ikik} - C_{iikk} 
  \right).
\end{equation}
The two second-order deviators can be calculated by 
\begin{align}
  &D_{ij} = \frac{5}{7} \left( C_{kkjj}-\frac{1}{3}   
    C_{kkll}\delta_{ij}\right) - \frac{4}{7} \left( 
    C_{kiki}-\frac{1}{3} C_{klkl}\delta_{ij}\right),\\
  &\hat{D}_{ij} = \frac{3}{7} \left(   
    C_{kikj}-\frac{1}{3}C_{klkl}\delta_{ij}\right) - 
    \frac{2}{7} \left( C_{kkij} -\frac{1}{3} 
    C_{kkll}\delta_{ij}\right)
\end{align}
This allows to compute the fourth-order deviator by 
  removing the other parts of the deviatoric
  decomposition
\begin{equation}
 \begin{split}
  D_{ijkl}& = C_{ijkl}-\left( \lambda \delta_{ij}\delta_{kl} +
    \mu (\delta_{ik}\delta_{jl} + \delta_{il}\delta_{jk}) + 
    \delta_{ij}D_{kl}  \right.\\
  & \left. + \delta_{kl} D_{ij} + \left( \delta_{ik}\hat{D}_{jl}
    + \delta_{il}\hat{D}_{jk} + \delta_{jk}\hat{D}_{il} + 
    \delta_{jl}\hat{D}_{ik}\right) \right)
 \end{split}
\end{equation}
  with the Kronecker delta $\delta_{ij}$.

As a side comment, one can also use this decomposition
  to calculate Young's modulus $E(\textbf{d})$ in a 
  specified direction $\textbf{d}$. 
Basically, the Young's modulus describes the stiffness upon uniaxial stretching in this 
  direction. 
  B\"ohlke and Br\"uggemann \cite{bohlke} calculated it by 
\begin{equation}
    \frac{1}{E(\textbf{d})} = \frac{1}{E^{RI}} +  6 
    \textbf{D} \cdot \textbf{d} \otimes \textbf{d} + 
    \textbf{d} \otimes \textbf{d} \cdot \mathcal{D}: 
    (\textbf{d} \otimes \textbf{d})
\end{equation}
  where $E^{RI} = 1/(2\mu + \lambda)$ is used to normalize the
  directional dependent quantities.
  
\subsubsection{Multipole Decomposition} \label{sec:MultipoleDecomposition}
Starting with a stiffness tensor given by its 21 coefficients
  in arbitrary Cartesian coordinates, the deviatoric decomposition 
  \eqref{deviatoricDec} and the multipole representation 
  \eqref{eq:multipole} of the stiffness tensor can be 
  used to calculate the position of the symmetry planes of any 
  anisotropic material.

To calculate the multipoles the equation \eqref{eq:poly} is used.
For the second order deviators $\mathbf{D}=\mathbf{D}^1$ and
  $\hat{\mathbf{D}}=\mathbf{D}^2$, we set $q=2$ and
\begin{equation}
    \begin{split}
        &a_{2,0}^k = -\sqrt{\frac{3}{2}} (D_{11}^k + D_{22}^k ), \quad a_{2,1}^k = D_{13}^k - iD_{23}^k\\
        &a_{2,2}^k = \frac{1}{2} (D_{11}^k - D_{22}^k) -iD_{12}^k
        \end{split}
\end{equation}
  where $k=1,2$.
For the fourth-order deviator $\mathcal{D}$, we use $q=4$ and
\begin{equation}
    \begin{split}
        & a_{4,0} = \sqrt{\frac{35}{8}} \left( D_{1111} + D_{2222} + 2D_{1122}\right)\\
        & a_{4,1} = \sqrt{\frac{7}{2}} (-D_{2213}-D_{1113} + i(D_{2223}+D_{1123}))\\
        & a_{4,2} = \frac{\sqrt{7}}{2} (D_{2222} -D_{1111} + 2i (D_{2212}+ D_{1112}))\\
        & a_{4,3} = \frac{1}{\sqrt{2}} (D_{1113}-3D_{2213}-i(3D_{1123}-D_{2223}))\\
        & a_{4,4} = \frac{1}{4} D_{1111} + \frac{1}{4} D_{2222} - \frac{3}{2} D_{1122} 
                    + i(D_{2212} - D_{1112})
    \end{split}
\end{equation}
We need the complex roots of these three polynomials.
In the second-order case, we can use the known formula for quadratic
  equations to solve for its four roots.
In the fourth-order case, we need to use a numerical method to find
  the eight roots. 
We use Laguerre's method to calculate the roots given by Press et al.
  \cite{press:1992}.
The multipoles can be calculated by using \eqref{eq:zeros} and
  \eqref{eq:mps}.

\subsubsection{Anisotropy Type}
One of the useful applications of the multipole decomposition is to compute the anisotropy
  type of the stiffness tensor.
For that purpose, the symmetry planes of each deviator must be calculated.
The intersection of these are the symmetry planes of the stiffness 
  tensor and determine the anisotropy type.
\begin{figure}[h]
    \centering
    \includegraphics[scale = 0.13]{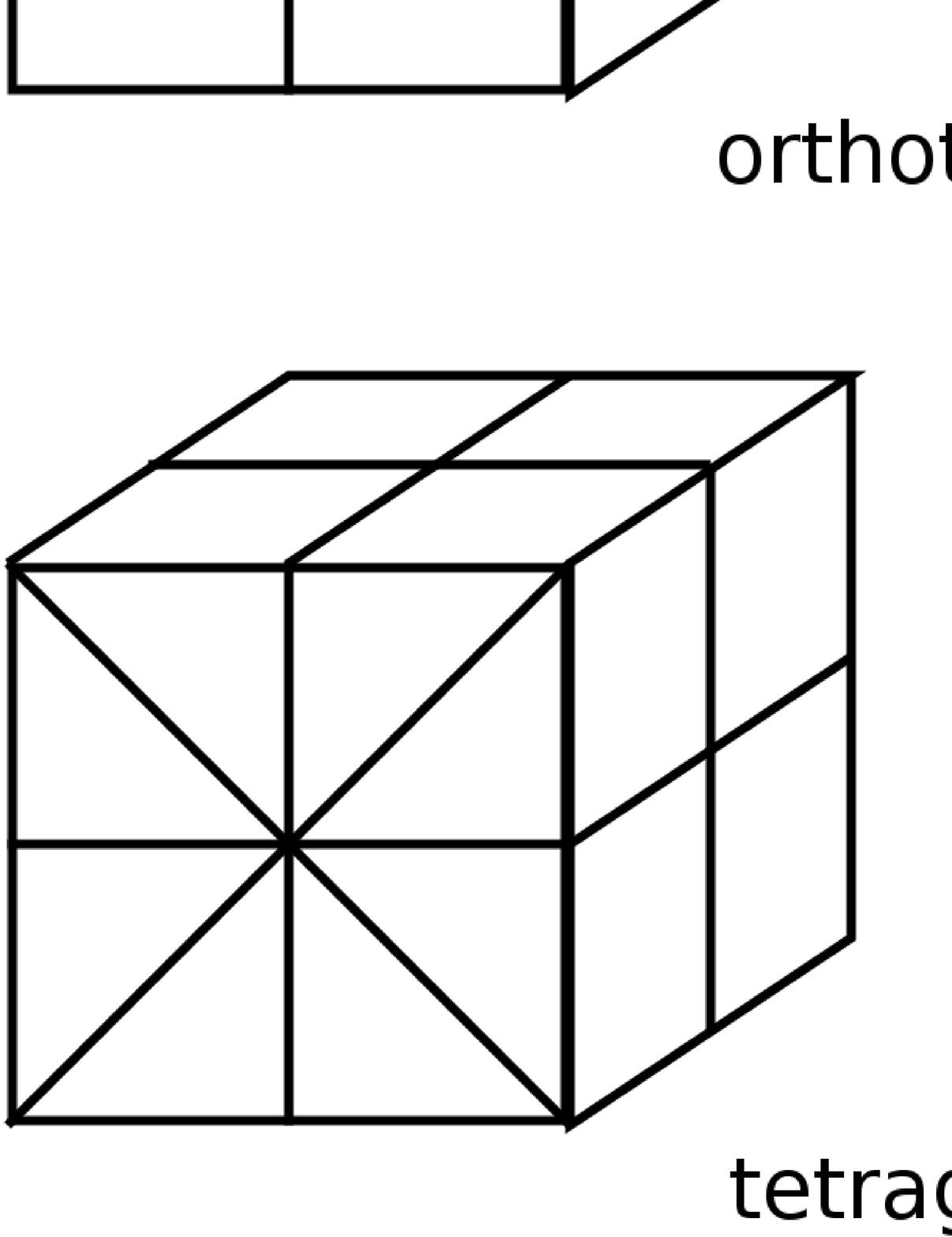}
    \caption{Anisotropy Classes of the stiffness tensor. Left: The symmetry planes of the material in this point. Right: The symmetries of the multipoles (yellow: fourth-order, blue: second-order).}
    \label{fig:aniso}
\end{figure}

A general stiffness tensor describes a rather complicated
  relation between stress and strain despite its linear form due to its incremental definition.
In practice, most engineering materials exhibit symmetries
  that simplify this relation further.
Of special interest are materials showing symmetry change under load.

In this context, a plane symmetry means that the elastic behavior,
  i.e. the stress-strain relation does not change under a 
  reflection at this plane.
The simplest materials are \textbf{isotropic}.
In this case, any plane is a symmetry plane and the relation between
  strain and stress is the same for all directions.
In this case, one needs only the two Lam\'{e} coefficients 
  in Eq.~\eqref{eq:lame} to describe the relation.
In a linear setting, $\lambda$ enters the scalar relation between 
  a uniform compression and isostatic pressure, and
$\mu$ describes the scalar relation between any volume-preserving
  (i.e. isochoric or deviatoric) strain and isochoric stress with
  the same direction.

In the general case, engineers distinguish the following 
  different classes of symmetries: \textbf{isotropic}, 
  \textbf{transversally isotropic}, \textbf{cubic}, \textbf{tetragonal}, 
  \textbf{orthotropic}, \textbf{monoclinic}, \textbf{trigonal} and \textbf{triclinic} materials. 
Of course, these characterizations hold only pointwise if the material
  exhibits different symmetry types at different positions.
In the following, we give a short characterization of all classes.
A more comprehensive treatment can be found in textbooks on solid 
  mechanics such as \cite{cowin2007}.
  
A second-order deviator can either be isotropic, transversally isotropic 
  or orthogonal. 
The set of symmetry plane normals $\text{MP}[\textbf{D}^1]$ in these
  three cases is given by
{\footnotesize
\begin{equation}
    \text{MP}[\textbf{D}^1] = \begin{cases}
        \text{all directions} & \text{if }A_1 =0\\
        \{ \textbf{n}, \textbf{m}_\varphi, \varphi \in [0,2\pi)\} & \text{if } 
        A_1 \neq 0,\, \textbf{n}_1 = \textbf{n}_2 = \textbf{n}\\
        \mathbf{N}_1, \mathbf{N}_2, \mathbf{N}_3 & \text{if }A_1\neq 0, 
        \,\textbf{n}_1 \neq \textbf{n}_2
    \end{cases}
\end{equation}}
where the infinite set $\{\textbf{m}_{\varphi}\}$ contains all vectors
  orthogonal to $\mathbf{n}$,
  and $\textbf{N}_1, \, \textbf{N}_2, \,\textbf{N}_3$ are given by 
\begin{equation}
    \textbf{N}_1 = \frac{\textbf{n}_1 + \textbf{n}_2}{\left| \textbf{n}_1 
    + \textbf{n}_2 \right|}, \quad \textbf{N}_2 = \frac{\textbf{n}_1 \times 
    \textbf{n}_2}{\left| \textbf{n}_1 \times \textbf{n}_2 \right|}, \quad 
    \textbf{N}_3 = \textbf{N}_1 \times \textbf{N}_2
\end{equation}
The set $\text{MP}[\textbf{D}^2]$ can be calculated equivalently using the
  multipoles $\textbf{n}_3,\textbf{n}_4$. 

The set $\text{MP}[\mathcal{D}]$ is a bit more complicated. 
We follow the description in Zou et al.~\cite{zou}.
If the deviator $\mathcal{D}$ is transversally isotropic the normals 
  are all the same $\textbf{n}_5 = \textbf{n}_6 = \textbf{n}_7 =
  \textbf{n}_8= \textbf{n}$. 
The symmetry plane normals can be calculated like these in the 
  second-order transversally isotropic deviator case.\\
If the deviator has cubic symmetry, the cube is given by
  a right-hand coordinate system
  $\left\lbrace \textbf{m}_1, \textbf{m}_2, \textbf{n} \right\rbrace$ 
  such that the axis direction set $\left\lbrace \textbf{n}_5, 
  \textbf{n}_6, \textbf{n}_7, \textbf{n}_8 \right\rbrace$ is given by
\begin{equation}
    \begin{split}
        \left\lbrace \frac{\textbf{n}+\textbf{m}_1 + \textbf{m}_2}{\sqrt{3}}, 
        \frac{\textbf{n}-\textbf{m}_1 + \textbf{m}_2}{\sqrt{3}}, \frac{\textbf{n}+ 
        \textbf{m}_1 - \textbf{m}_2}{\sqrt{3}}, \frac{-\textbf{n}+ \textbf{m}_1 + 
        \textbf{m}_2}{\sqrt{3}} \right\rbrace 
    \end{split}
\end{equation}
  and the nine symmetry planes are the three coordinate
  planes and the six planes created by rotating in each of the
  coordinate planes by an angle of $\pm \pi/4$.
If the deviator is tetragonal, the four multipoles are the 
  result of a rotation at $90^\circ$ around some axis.
The symmetry plane normals of $\mathcal{D}$ are given by 
  the vector orthogonal to all multipoles, the multipoles 
  themselves and the vectors that lie with an angle of $45^\circ$ 
  between two of the multipoles.
If the deviator is trigonal, one normal is orthogonal to 
  the other three, while the other three are coplanar.
The coplanar ones result from a rotation around the fourth multipole
  at $120^\circ$.
The three coplanar multipoles are the symmetry plane normals.
If the deviator is orthogonal, the multipoles result from a 
  rotation of one of the following sets
\begin{equation}
    \begin{cases}
        S_W =& \{ \textbf{n}(\theta, \varphi), \textbf{n}(\theta, \pi - \varphi), \textbf{n}(\theta, \pi + \varphi), \textbf{n}(\theta, 2\pi - \varphi)\};\\
        S_U =& \{ \textbf{n}(\pi /2, \varphi), \textbf{n} (\pi /2, \pi -\varphi) \} \text{ or } \{ \textbf{n}(\theta, 0), \textbf{n}(\theta, \pi)\} \\
        & \text{ or } \{ \textbf{n} (\theta, \pi /2) , \textbf{n}(\theta, 3\pi /2 )\};\\
        S_V =& \{ \text{e}_1\} \text{ or } \{ \textbf{e}_2 \} \text{ or } \{\textbf{e}_3 \}
    \end{cases}
\end{equation}
  and the symmetry plane normals equal the coordinate surface normals.
If all multipoles lie on a plane or there exists 
  a plane that is their mid-separating surface, the normal to this 
  plane is the symmetry plane normal. 
A triclinic deviator has four arbitrary multipoles not fulfilling 
  any of the above relations and there is no symmetry plane.
Finally, the symmetry plane normals of the stiffness 
  tensor are given by the intersection of the symmetry plane 
  normals of three deviators.
\begin{equation}\label{MP}
    \text{MP}[\mathcal{C}] = \text{MP}[\textbf{D}] \cap \text{MP}
    [\hat{\textbf{D}}] \cap \text{MP}[\mathcal{D}]
\end{equation} 

\section{Conclusion}
Every tensor of arbitrary order up to dimension three can be described
  uniquely by a set of vectors and scalars.
This is the main statement of this work.
This decomposition is not new, but it is neither well known nor 
  well understood.
Thus, this general concept has not found its way into natural or engineering sciences and 
its physical interpretation remains unclear in many areas.
As a step towards obtaining such interpretations, we gave a connection between the multipole decomposition and the
  spectral decomposition of a symmetric three-dimensional
  second-order tensor.
We are sure, there are other interesting facts about this 
  decomposition and we want to motivate more researchers to 
  analyze this decomposition and use it for different examples to 
  explore its application-specific meanings.

\bibliographystyle{unsrt}  

\bibliography{egbibsample}

\begin{thebibliography}{10}

\bibitem{laidlaw2012new}
David~H Laidlaw and Anna Vilanova.
\newblock {\em New developments in the visualization and processing of tensor
  fields}.
\newblock Springer Science \& Business Media, 2012.

\bibitem{hotz2015visualization}
Ingrid Hotz and Thomas Schultz.
\newblock {\em Visualization and processing of higher order descriptors for
  multi-valued data}.
\newblock Springer, 2015.

\bibitem{kratz2013visualization}
Andrea Kratz, Cornelia Auer, Markus Stommel, and Ingrid Hotz.
\newblock Visualization and analysis of second-order tensors: Moving beyond the
  symmetric positive-definite case.
\newblock In {\em Computer Graphics Forum}, volume~32, pages 49--74. Wiley
  Online Library, 2013.

\bibitem{itskov2007tensor}
Mikhail Itskov.
\newblock {\em Tensor algebra and tensor analysis for engineers}.
\newblock Springer, 2007.

\bibitem{Maxwell:1881}
James~Clerk Maxwell.
\newblock {\em A treatise on electricity and magnetism}, volume~1.
\newblock Clarendon press, 1881.

\bibitem{backus}
George Backus.
\newblock A geometrical picture of anisotropic elastic tensors.
\newblock {\em Reviews of geophysics}, 8(3):633--671, 1970.

\bibitem{zou2003}
W-N Zou and Q-S Zheng.
\newblock Maxwell's multipole representation of traceless symmetric tensors and
  its application to functions of high-order tensors.
\newblock In {\em Proceedings of the Royal Society of London A: Mathematical,
  Physical and Engineering Sciences}, volume 459, pages 527--538. The Royal
  Society, 2003.

\bibitem{zou2001}
W-N Zou, Q-S Zheng, D-X Du, and J~Rychlewski.
\newblock Orthogonal irreducible decompositions of tensors of high orders.
\newblock {\em Mathematics and Mechanics of Solids}, 6(3):249--267, 2001.

\bibitem{Kelvin:1856}
W~Thomsen~Lord Kelvin.
\newblock Elements of a mathematical theory of elasticity, part 1: On stresses
  and strains.
\newblock {\em Philosophical Transactions of the Royal Society}, 166:481--498,
  1856.

\bibitem{Mehrabadi1990}
Morteza~M. Mehrabadi and Stephen~C. Cowin.
\newblock {Eigentensors of linear anisotropic elastic materials}.
\newblock {\em The Quarterly Journal of Mechanics and Applied Mathematics},
  43(1):15--41, 02 1990.

\bibitem{Nagel2016}
Thomas Nagel, Uwe-Jens G{\"{o}}rke, Kevin~M. Moerman, and Olaf Kolditz.
\newblock {On advantages of the Kelvin mapping in finite element
  implementations of deformation processes}.
\newblock {\em Environmental Earth Sciences}, 75(11):937, jun 2016.

\bibitem{wigner}
E.P. Wigner.
\newblock Group theory.
\newblock {\em Academic Press, New York}, 1959.

\bibitem{Sylvester:1876}
James~Joseph Sylvester.
\newblock Note on spherical harmonics.
\newblock {\em The London, Edinburgh, and Dublin Philosophical Magazine and
  Journal of Science}, 2(11):291--307, 1876.

\bibitem{zou}
W-N Zou, C-X Tang, and W-H Lee.
\newblock Identification of symmetry type of linear elastic stiffness tensor in
  an arbitrarily orientated coordinate system.
\newblock {\em International Journal of Solids and Structures},
  50(14-15):2457--2467, 2013.

\bibitem{Zou:2003}
W-N Zou and Q-S Zheng.
\newblock Maxwell's multipole representation of traceless symmetric tensors and
  its application to functions of high-order tensors.
\newblock {\em Proceedings: Mathematics, Physical and Engineering Sciences},
  pages 527--538, 2003.

\bibitem{zheng2006degenerate}
Xiaoqiang Zheng, Xavier Tricoche, and Alex Pang.
\newblock Degenerate 3d tensors.
\newblock In {\em Visualization and Processing of Tensor Fields}, pages
  241--256. Springer, 2006.

\bibitem{hergl2019visualization}
Chiara Hergl, Thomas Nagel, Olaf Kolditz, and Gerik Scheuermann.
\newblock Visualization of symmetries in fourth-order stiffness tensors.
\newblock In {\em 2019 IEEE Visualization Conference (VIS)}, pages 291--295.
  IEEE, 2019.

\bibitem{Weyl:1946}
Hermann Weyl.
\newblock {\em The classical groups}.
\newblock Princeton University Press, 1946.

\bibitem{bohlke}
Thomas B{\"o}hlke and C~Br{\"u}ggemann.
\newblock Graphical representation of the generalized hooke's law.
\newblock {\em Technische Mechanik}, 21(2):145--158, 2001.

\bibitem{press:1992}
William~H Press, Saul~A. Teukolsky, William~T. Vetterling, and Brian~P.
  Flannery.
\newblock Numerical recipes in c: The art of scientific computing.
\newblock 1992.

\bibitem{cowin2007}
Stephen~C Cowin and Stephen~B Doty.
\newblock {\em Tissue mechanics}.
\newblock Springer Science \& Business Media, 2007.

\end{thebibliography}

\end{document}